\documentclass[12]{amsart}
\usepackage[utf8]{inputenc}
\usepackage[a4paper, total={6in, 8.5in}]{geometry}
\usepackage{hyperref}
\usepackage{url}
\providecommand{\myfloor}[1]{ \lfloor #1 \rfloor }
\newtheorem{theorem}{Theorem}

\newcommand{\QEDB}{\hfill\ensuremath{\square}}
\title{A Note on the Number of Regions in a Line Arrangement}
\author{Dickson Y. B. Annor, Michael S. Payne }

\address{Department of Mathematics and Physics, La Trobe University, Bendigo, Victoria 3552 Australia.}
\email{d.annor@latrobe.edu.au, m.payne@latrobe.edu.au}
\thanks{Dickson Annor is supported by La Trobe Graduate Research Scholarship. Michael Payne is supported by a DECRA from the Australian Research Council.}

\begin{document}
\maketitle

\maketitle
\begin{abstract}
For an  arrangement of $n$ lines in the real projective plane, we denote by $f$ the number of  regions into which the real projective plane is divided by the lines. Using Bojanowski's inequality, we establish a new lower bound for $f$. In particular, we show that if no more than $\frac{2}{3}n$ lines intersect at any point, then $f \ge \frac{1}{6}n^{2}$. 
\end{abstract}

\vspace{0.5cm}

\maketitle

\maketitle




\section{Introduction}\label{sec1}
 Let $\mathcal{L}$ be an arrangement of $n \ge 2$ lines in the real projective plane  $\mathbb{RP}^{2}$ and let $m$ denote  the maximum number of lines from $\mathcal{L}$ intersecting at one point. The lines from $\mathcal{L}$ divide $\mathbb{RP}^{2}$ into polygonal regions which are the connected components of the complement of the union of the lines. Denote the number of regions by $f$. The question we are interested in is: how many regions can be obtained (under all possible arrangements $\mathcal{L}$ of $n$ lines)?\\
 

Below we collect some known lower bounds for $f$ in terms of $n$  and $m$.

\begin{itemize}
  \item $f \ge 2n - 2$, if $m < n$, \hfill  Gr\"{u}nbaum \cite{MR0307027};
  \item $f \ge 3n - 6$, if $m \le n - 2$, \hfill Gr\"{u}nbaum \cite{MR0307027};
\item $f \ge m(n + 1 - m)$, \hfill Arnol'd \cite{arnol2008into};
\item  $f \ge \frac{n(n - 1)}{2(m -1)}$, if $m > 2$, \hfill Arnol'd \cite{arnol2008into};
\item  $f \ge (m + 1)(n - m)$, \hfill Arnol'd \cite{arnol2008into} and  Purdy \cite{MR566442};
\item  $f \ge (r +1)(n - r)$, if  $m \le n - r$ and $n \ge \frac{r^{2} + r}{2} + 3$ for some $r \in \mathbb{Z}$,
\hfill Shnurnikov \cite{MR2815265};
\item $f \ge  2\left( \frac{n^{2}- n + 2m}{m + 3}\right)$, \hfill  Shnurnikov \cite{MR2815265};
\item  $f \ge \frac{(3m -10)n^{2} +(m^{2} -6m +12)n}{m^{2} + 3m -18} + 1$, if $5 \le m < n - 2$, \hfill  Shnurnikov \cite{MR3458599}.\\
\end{itemize}
 
 In this paper we use  Bojanowski's inequality \cite{bojanowski2003zastosowania} to  establish a new lower bound for $f$.  The main result states that:
 
 \begin{theorem}\label{t1}
Let $\mathcal{L}$ be an arrangement of $n$ lines in the real projective plane such that $m \le \frac{2}{3}n$. Then 
$$
f \ge \frac{(m + 2)n^{2} + (3m - 6)n}{6m} + 1 \ge \frac{1}{6}n^{2}.
$$ 
\end{theorem}

We remark that to the best of our knowledge, if $m(n)$ is a sublinear but increasing function of $n$, then this is the first quadratic lower bound on $f$. For example consider the case $m = \sqrt{n}$ in the previously known inequalities given above.


 
  \section{Bounds for Number of Regions}
  
  For an arrangement of lines $\mathcal{L}$ in the projective plane we denote by  $t_{k}$, $2 \le k \le n- 1$,  the number of intersection points  where exaclty $k$ lines of the arrangement are incident. The following are some known relations for values of $t_{k}$.

\begin{itemize}
  \item $ t_{2} \ge 3 +  \sum_{k \ge 4} (k - 3) t_{k}$, \hfill Melchior \cite{MR4476};
  \item $ t_{2} \ge \frac{6}{13}n$ for $n \ge 8$, \hfill Csima and Sawyer \cite{MR1194036};
\item $ t_{2} + \frac{3}{4}t_{3} \ge n +  \sum_{k \ge 5} (2k - 9) t_{k}$, if $t_{n - 1} = t_{n - 2} = 0$,\hfill Hirzebruch \cite{MR860410};
\item $ t_{2} + \frac{3}{4}t_{3} \ge n +  \sum_{k \ge 5} \left(\frac{1}{4}k^{2} - k\right) t_{k}$, if $t_{k} = 0$ for $k > \frac{2}{3}n$,\hfill   Bojanowski \cite{bojanowski2003zastosowania};
\item $t_{2} \ge \frac{1}{2}n$ and $t_{2} \ge 3\myfloor{\frac{1}{4}n}$ for sufficiently large, even and odd $n$, respectively, Green and Tao \cite{MR3090525}.
\end{itemize}
Maybe it is worth to mention here that both  Bojanowski  \cite{bojanowski2003zastosowania} and Hirzebruch \cite{MR860410} inequalities hold for arrangements of complex lines in the complex projective plane and consequently, they also hold for arrangements of lines in the real projective plane. To the best of our knowledge,  Bojanowski's inequality \cite{bojanowski2003zastosowania} is the strongest known inequality for line arrangements with $m \le  \frac{2}{3}n$ \cite{MR4245242}.
  \\\\
  \textit{Proof of Theorem  \ref{t1}.}
  Let  $\mathcal{L}$ be an arrangement of $n$ lines. If we add lines one by one, then the number of new regions created by each line is equal to the number of intersection points with previously added lines. In this process, a point with $k$ lines passing through it is intersected $k - 1$ times. Thus, the number of regions, including 1 for the first line, is 
\begin{equation}\label{n1}
\begin{split}
f = 1 +  \sum_{k = 2}^{m}(k - 1) t_{k}.
\end{split}
\end{equation}
Note that (\ref{n1}) can be obtained by using the fact that the Euler characteristic of the real  projective plane is 1. The number of pairs of lines in  $\mathcal{L}$ is equal to $\frac{n(n-1)}{2}$. In a projective plane, every pair of lines intersects at exactly one point, and if $k$ lines meet at a point, we get  $\frac{k(k-1)}{2}$ of such pairs which cross at that point. Since $t_{k} = 0$ for $k > m$, we obtain
\begin{equation}\label{n2}
n(n-1) = \sum^{m}_{k = 2}k(k - 1) t_{k}.
\end{equation}
Suppose we are given an inequality
\begin{equation}\label{r3}
 \sum_{k = 2}^{m}\alpha_{k} t_{k} \ge \alpha_{0}
\end{equation}
where $\alpha_{0}, \alpha_{2}, \alpha_{3}, \cdots, \alpha_{m}$ are some  real numbers, and suppose that for some $c_{1}, c_{2} > 0$ the inequality

\begin{equation} \label{r1}
 c_{1}k(k-1) + c_{2}\alpha_{k} \le k - 1
\end{equation}
is satisfied for all $2 \le k \le m$. Multiply both sides of (\ref{r1}) by $t_{k}$ and sum up for $k = 2, 3, \cdots , m$ to obtain
\begin{equation*}
c_{1} \sum^{m}_{k = 2}k(k - 1) t_{k} + c_{2} \sum^{m}_{k = 2} \alpha_{k} t_{k} \le   \sum^{m}_{k = 2}(k - 1) t_{k}
\end{equation*}
since $t_{k} \ge 0$.
This  is equivalent to 
\begin{equation}\label{r2}
c_{1}n(n-1) + c_{2}\sum^{m}_{k = 2} \alpha_{k} t_{k} \le  f - 1 .
\end{equation}
Using (\ref{r3}) and (\ref{r2}) and the fact that $c_{2} > 0$, we obtain
\begin{equation}\label{r4}
f \ge c_{1}n(n - 1) + c_{2}\alpha_{0} + 1
\end{equation}
for positive $c_{1}, c_{2},$ satisfying (\ref{r1}).
For $m \le \frac{2}{3}n$ we use Bojanowski's inequality \cite{bojanowski2003zastosowania} 
\begin{equation*}
t_{2}  + \frac{3}{4}t_{3} + \sum_{k \ge 5} \left(k - \frac{1}{4}k^{2}\right) t_{k}  \ge n 
\end{equation*}
in the form (\ref{r3}) to obtain the following
\begin{equation*}
\alpha_{0} = n,\;\;\ \alpha_{2} = 1, \;\; \alpha_{3} = \frac{3}{4},\;\; \alpha_{4} = 0 ,\;\; \alpha_{k} = \left(k- \frac{1}{4}k^{2}\right) \;\; \mbox{for} \; k \ge 5.
\end{equation*}
From (\ref{r1}) we get 
\begin{align*}
& 1 \geq 2c_{1} + c_{2}, \;\;\;\; 2 \geq 6c_{1} + \frac{3}{4}c_{2}, \;\;\;\; 3 \geq 12c_{1},\;\;\mbox{for} \;\; k = 2, 3, \mbox{and} \; 4, \mbox{respectively},\\
& 0 \ge c_{1} k(k - 1) +  c_{2}\left(k - \frac{1}{4}k^{2}\right) - (k - 1) \;\;\; \mbox{for} \;\;\; 5 \le k \le m.
\end{align*}
For $m \ge 2$, let us take the positive numbers
\begin{equation*}
 c_{1} = \frac{m + 2}{6m }, \;\;\;\;\;\;\;\;\;\;\;\;\;\;  c_{2} = \frac{2(m - 1)}{3m}.
\end{equation*}
Now we need to check these inequalities for $2 \le k \le m$ and for the given $c_{1}, c_{2}$. The first three are easy to check, so we verify the last one for $5 \le k \le m$. Thus,
\begin{equation*}
  c_{1} k(k - 1) +  c_{2}\left(k -\frac{1}{4}k^{2}\right) - (k - 1) = \dfrac{0.5(k - 2)(k - m)}{m} \le 0,
\end{equation*}
because $k \le m$ and $k \ge 5$. So, we obtain (\ref{r4}) for the given $c_{1}, c_{2},$ and hence, the inequality of the theorem.\QEDB\\

Note that lower bounds on $f$ in the form (\ref{r4}) were obtained by Shnurnikov in \cite{MR3458599}.  In \cite{MR3458599} he applied Hirzebruch's inequality \cite{MR860410} to obtain the result mentioned in Section  \ref{sec1}.\\


It is natural to ask under which assumptions the inequality of Theorem  \ref{t1} is stronger than previously known inequalities. The inequality in Theorem  \ref{t1} is quadratic in $n$. So, it sufficies to compare it to those inequalities mentioned in Section  \ref{sec1} that are quadratic in $n$ for some function $m(n)$. In particular, the results of Arnol'd \cite{arnol2008into} and Purdy \cite{MR566442} become quadratic in $n$ if $m(n) = \frac{n}{p}$ where $p$ is a real number greater than 1. A simple calculation shows that Theorem  \ref{t1} is weaker than those inequalities  when $p \in \left(3 - \sqrt 3, 3 + \sqrt 3\right)$ and it is also weaker than  Shnurnikov \cite{MR2815265}  when $m  \le 5$. On the other hand, Theorem  \ref{t1} is stronger than all the inequalities mentioned in Section \ref{sec1} whenever $7 \le m \le \dfrac{n}{5}$ and for $m = 6$ we have equality with Shnurnikov \cite{MR3458599}. 

\section*{Acknowledgements}
We are grateful to Yuri Nikolayevsky for his interest in this work and for his useful comments and suggestions which have improved the presentation. 

 
 
 \bibliographystyle{plain}
 \bibliography{bibliography.bib}

\begin{thebibliography}{10}

\bibitem{arnol2008into}
V.~I. Arnol’d.
\newblock Into how many regions do $n$ lines divide the plane?
\newblock {\em Matem. Pros., Ser. 3}, 12:95, 2008.

\bibitem{bojanowski2003zastosowania}
R.~Bojanowski.
\newblock {\em Zastosowania uog{\'o}lnionej nier{\'o}wnosci
  Bogomolova-Miyaoka-Yau}.
\newblock PhD thesis, Master Thesis (in Polish),
  {\url{http://www.mimuw.edu.pl/\%7Ealan/postscript}}, 2003.

\bibitem{MR1194036}
J.~Csima and E.~T. Sawyer.
\newblock There exist {$6n/13$} ordinary points.
\newblock {\em Discrete Comput. Geom.}, 9(2):187--202, 1993.

\bibitem{MR3090525}
B.~Green and T.~Tao.
\newblock On sets defining few ordinary lines.
\newblock {\em Discrete Comput. Geom.}, 50(2):409--468, 2013.

\bibitem{MR0307027}
B.~Gr\"{u}nbaum.
\newblock {\em Arrangements and spreads}.
\newblock Conference Board of the Mathematical Sciences Regional Conference
  Series in Mathematics, No. 10. American Mathematical Society, Providence,
  R.I., 1972.

\bibitem{MR860410}
F.~Hirzebruch.
\newblock Singularities of algebraic surfaces and characteristic numbers.
\newblock In {\em The {L}efschetz centennial conference, {P}art {I} ({M}exico
  {C}ity, 1984)}, volume~58 of {\em Contemp. Math.}, pages 141--155. Amer.
  Math. Soc., Providence, RI, 1986.

\bibitem{MR4476}
E.~Melchior.
\newblock \"{U}ber {V}ielseite der projektiven {E}bene.
\newblock {\em Deutsche Math.}, 5:461--475, 1941.

\bibitem{MR4245242}
P.~Pokora.
\newblock Hirzebruch-type inequalities viewed as tools in combinatorics.
\newblock {\em Electron. J. Combin.}, 28(1):Paper No. 1.9, 22, 2021.

\bibitem{MR566442}
G.~Purdy.
\newblock On the number of regions determined by {$n$} lines in the projective
  plane.
\newblock {\em Geom. Dedicata}, 9(1):107--109, 1980.

\bibitem{MR2815265}
I.~N. Shnurnikov.
\newblock Into how many regions do {$n$} lines divide a plane if at most
  {$n-k$} of them are collinear?
\newblock {\em Vestnik Moskov. Univ. Ser. I Mat. Mekh.}, (5):32--36, 2010.

\bibitem{MR3458599}
I.~N. Shnurnikov.
\newblock A {$t_k$} inequality for arrangements of pseudolines.
\newblock {\em Discrete Comput. Geom.}, 55(2):284--295, 2016.

\end{thebibliography}
\end{document}